\documentclass[final,12pt]{elsarticle}

\usepackage{amsmath, amsthm,amsfonts}
\usepackage{tikz}

\usepackage{float}
\restylefloat{figure}

\newcommand{\A}{T}
\newcommand{\B}{V}
\newcommand{\C}{\mathcal{C}}
\newcommand{\landO}{\mathcal{O}}
\newtheorem{theor}{Theorem}

\DeclareMathOperator{\e}{e}
\newcommand{\dd}{{\rm d}}
\journal{Applied Numerical Mathematics}

\begin{document}
\begin{frontmatter}
\title{Force-Gradient Nested Multirate Methods for Hamiltonian Systems}

\author{Dmitry Shcherbakov\corref{cor1}} 
\ead{shcherbakov@math.uni-wuppertal.de}
\cortext[cor1]{corresponding author}

\author{ Matthias Ehrhardt}
\ead{ehrhardt@math.uni-wuppertal.de}

\author{Michael G\"unther}
\ead{guenther@math.uni-wuppertal.de}

\address{Lehrstuhl Angewandte Mathematik und Numerische Analysis, Bergische Universit\"{a}t Wuppertal,
Gau{\ss}strasse 20, 42119 Wuppertal, Germany}

\author{Michael Peardon}
\ead{mjp@maths.tcd.ie}
\address{School of Mathematics, Trinity College, Dublin 2, Ireland}

\begin{abstract}
Force-gradient decomposition methods are used to improve the energy preservation of 
symplectic schemes applied to Hamiltonian systems. 
If the potential is composed of different parts with strongly varying dynamics, this multirate potential can 
be exploited by coupling force-gradient decomposition methods
with splitting techniques for multi-time scale problems to further increase the accuracy of the scheme and reduce 
the computational costs. 
In this paper, we derive novel force-gradient nested methods and test them numerically. Such methods can be used to increase 
the  acceptance rate for the molecular dynamics step of the Hybrid Monte Carlo 
algorithm (HMC)
and hence improve its computational efficiency.
\end{abstract}


\end{frontmatter}

\section{Introduction}
For classical mechanical systems, the equation of motion can be written as
\begin{equation}\label{eq:class_em}
\frac{\dd\mathbf{\rho}}{\dd t} =[\rho \circ H ]\equiv \mathcal{L}(t) \mathbf{\rho}(t),
\end{equation}
where $\mathbf{\rho}$  is the set of phase variables, $[~  \circ ~]$ denotes the Poisson bracket, $H$
 represents the Hamiltonian function, and $\mathcal{L}$ denotes the Liouville operator. 
For the case of $N$ particles, located in a spatially inhomogeneous
 time-dependent external field $u(\mathbf{r}_i,t)$ and interacting through 
the pair-wise potential 
$\varphi(r_{ij})\equiv \varphi (| \mathbf{r}_i-\mathbf{r}_j|)$, the {\em Hamiltonian} reads
\begin{equation}\label{eq:hamiltonian}
 H= \sum\limits^{N}_{i=1} \frac{m_i \mathbf{v}^{2}_{i}}{2} +\frac{1}{2}\sum\limits^{N}_{i\neq j} \varphi(r_{ij})
+  \sum\limits^{N}_{i=1}u(\mathbf{r}_i,~t) \equiv T(\mathbf{v}) + V(\mathbf{r}).
\end{equation}
Here $\mathbf{r}_i$ represents the position of particle $i$ $(i=1,2,\dots,N)$ moving with velocity 
$\mathbf{v}_i=\dd \mathbf{r}_i/\dd t $ and carrying the mass $m_i$, so that $T$ and $V$ are the total kinetic and
 potential energies, respectively. 
Then $\mathbf{\rho} =\{\mathbf{r}_i,~\mathbf{v}_i\} 
\equiv \{\mathbf{r},~\mathbf{v}\}$, 
and the {\em Liouville operator} of the system takes the form
\begin{equation}\label{eq:liouville}
 \mathcal{L}(t) = \sum\limits^{N}_{i=1} \left( \mathbf{v}_i \cdot \frac{\partial}{\partial \mathbf{r}_i}
 + \frac{\mathbf{f}_i(t)}{m_i} \cdot \frac{\partial}{\partial \mathbf{v}_i}  \right),
\end{equation}
where 
\begin{equation*}
\mathbf{f}_i(t) = \sum^{N}_{j(j\neq i)} \frac{\varphi '(r_{ij}) \mathbf{r_{ij}} }{r_{ij}} 
- \frac{\partial u(\mathbf{r}_i,t)}{\partial \mathbf{r}_i}
\end{equation*}
 are forces acting on the particles due to their interactions.

If the initial configuration $\mathbf{\rho}(0)$ is specified, the unique solution to the problem of Eqn.~\eqref{eq:class_em} 
can be presented by the {\em time propagator operator} as
\begin{equation}\label{eq:unique_solution}
 \rho(t)  = \left[ \e^{(\mathcal{D} + \mathcal{L}) h } \right]^{l} \rho(0), 
\end{equation}
where $h$ is a temporal step size and $l= t/h$ the total number of steps.
 $\mathcal{D} = \overleftarrow{\partial} / \partial t$ denotes the time derivative operator
acting on the left of  time-dependent functions. 
If $\mathcal{L}$ does not depend explicitly on time we set $\mathcal{D}=0$. 
In case of many-particle systems $(N > 2)$ the time propagator cannot be computed 
exactly even in the absence of time dependent potentials. 
Hence one has to apply numerical integration methods such as decomposition schemes, which both preserve the physical
 properties of the Hamiltonian system \eqref{eq:unique_solution} 
(symplecticity, time reversibility) and are computationally  efficient \cite{OmMrFo03}. 

The basic idea of a decomposition approach is to factor out the exponential propagator 
 $e^{(\mathcal{D}+\mathcal{L})h}$ in \eqref{eq:unique_solution}, 
such that $\mathcal{D} + \mathcal{L} = \hat{\A} + \hat{\B}$, 
where the differential operators $\hat{\A}=\mathbf{v} \cdot \partial /\partial\mathbf{r}$
and $\hat{\B}=\mathbf{a}\cdot \partial/\partial\mathbf{v}$ represent the kinetic and potential energies 
 with the acceleration  $\mathbf{a}=\{ \mathbf{a}_{i} \} =\{ \mathbf{f}_{i} / m_{i}\}$. 

In Section 2 we will discuss two well-known approaches to increase the accuracy of decomposition schemes:
 force-gradient schemes and nested multirate algorithms, which are both based on decomposition techniques. 
For both methods, computing the shadow Hamiltonian is the suitable tool for deriving the order of the numerical integration scheme. 
In Section 3 both approaches are combined to obtain a more efficient scheme. Finally, 
numerical results for a three body problem confirm the theoretical findings in Section 4.


\section{Methods for obtaining Higher Order Schemes}
In this section we will briefly recapitulate two well-known schemes
 (force-gradient and nested multirate schemes) for reducing the computational costs. 
As both approaches are based on decomposition, the computation of the shadow Hamiltonian can be used to determine the order of the numerical scheme.

\subsection{Shadow Hamiltonians}
When transferring the well-known concept of modified equations to Hamiltonian systems one ends up with the Hamiltonian if 
and only if the integrator is symplectic \cite{SH01}. 
The motivation for studying numerically the conservation properties of these 'modified Hamiltonians' are multifaceted \cite{SH01},
 e.g.\
numerical evidence for the existence of a Hamiltonian for a particular calculation, 
exposure of energy drifts caused by numerical instability, etc..
Skeel and Hardy \cite{SH01} proposed a simple strategy 
for deriving highly accurate estimates for modified Hamiltonians. 
Since these modified Hamiltonians approximate well the true Hamiltonian,  they are referred as "shadow" Hamiltonians $\tilde{H}$, cf.\ \cite{GS00}.
 The existence of these shadow Hamiltonians guarantees
the boundedness of the error in the symplectic map, in fact we
have $\tilde{H}(\mathbf{r},\mathbf{v},h)\to H(\mathbf{r},\mathbf{v})$ for $h\to0$.

Conversely, if one starts from a given numerical solver then it is well known 
that any symplectic integrator different from the Hamiltonian flow 
itself does not preserve the Hamiltonian however a nearby system, 
the so-called \textit{shadow Hamiltonian} $\tilde{H}$ is conserved. The energy
computed from the shadow Hamiltonian of a symplectic integrators differs by  
$H(\mathbf{r},\mathbf{v}) - \tilde{H}(\mathbf{r},\mathbf{v},h) = \landO(h^{p})$ from the true 
Hamiltonian \cite{KeCl07}, with $p$ being the order of the integration scheme.
Hence, computing the shadow Hamiltonian of a symplectic integrator is 
equivalent to determining the order of the integrator.

To compute a shadow Hamiltonian it is necessary to expand an exponential map to a \textit{Hausdorff series}. 
To do so, we need to use the {\em Baker-Cambell-Hausdorf (BCH) formula} \cite{KeSiCl12}.
\begin{equation}\label{eq:BCH}
 \ln(\e^{\hat{\A}}\e^{\hat{\B}}) = \sum\limits^{\infty}_{n=1} c_{n}(\A,\B),
\end{equation}
where the coefficients $c_n$ are recursively determined from the relations $c_1 = \A+\B$ and
\begin{multline*}
 (n+1)c_{n+1} =  \\
                \sum\limits^{\lfloor n/2\rfloor}_{m=1} \frac{B_{2m}}{(2m)!} \sum\limits_{k_1,\ldots,k_{2m}\geq 1} \text{ad} c_{k_1} \ldots \text{ad} c_{k_{2m}}(\A+\B) - \frac{1}{2}
                ( \text{ad} c_{n} )(\A-\B),
\end{multline*}
for $n \geq 0$, where  $\text{ad}a :b \mapsto [a,b]$ and $B_n$ denote the Bernoulli numbers.
For example, the shadow Hamiltonian of the leap-frog method $\e^{h\frac{\A}{2}}\e^{h\B}\e^{h\frac{\A}{2}} $ is given by 
 \begin{equation*}
 \tilde{H} = H -\frac{h^{2}}{24} \Bigl(2\bigl[\B,[\A,\B]\bigr]+\bigl[\A,[\A,\B]\bigr]\Bigr) +\landO(h^{4}),
 \end{equation*}
 which is of second order accuracy.

\subsection{Force-gradient schemes}
Force-gradient schemes are based on the fact that 
the total propagator in Eqn.~\eqref{eq:unique_solution} can be split in the following way:
\begin{equation} \label{eq:decomposition}
 \e^{(\hat{\A}+\hat{\B})h+\landO(h^{K+1})} =\prod\limits_{p=1}^{P} \e^{\hat{\A} a_{p}h} \e^{\hat{\B} b_{p}h +\C c_{p}h^3},
\end{equation}
where $\C = [\B,[\A,\B]]$ and $[~,~]$ denotes the commutator of two operators.

The coefficients $a_{p}$, $b_{p}$ and $c_{p}$ in \eqref{eq:decomposition} have to be chosen 
in such way to obtain the highest possible  order $K\ge1$ for a given integer $P\ge1$.
Eqn.~\eqref{eq:decomposition} represents the general form of the decomposition,
while  for $c_p \equiv 0 $ the decomposition reduces to the standard  
non-gradient factorization.
The force-gradient method is defined by using the value of $c_{p}$ 
which reduces the difference between the true Hamiltonian and shadow 
Hamiltonian $\tilde{H}$ which is conserved by the method. 
We will show how to determine the shadow Hamiltonian $\tilde{H}$ in the next section. 

The third order force-gradient operator $\C$ can be obtained for classical 
systems and is given by
\begin{equation*}
 \C \equiv \bigl[\B,[\A,\B]\bigr] = \sum \limits^{N}_{i=1} \frac{\mathbf{g}_{i}}{m_{i}} \cdot \frac{\partial}{\partial \mathbf{v}} \equiv  \mathbf{G} \cdot \frac{\partial}{\partial \mathbf{v}},
\end{equation*}
where 
\begin{equation*}
 \mathbf{g}_{i\alpha} =2 \sum_{j\beta} \frac{\mathbf{f}_{j\beta}}{m_j} \frac{\partial \mathbf{f}_{i\alpha}}{ \partial\mathbf{r}_{j\beta}} , 
\end{equation*}
$\alpha$ and $\beta$ denote the Cartesian components of the vectors.  
The force-gradient evaluations 
$\partial \mathbf{f}_{i\alpha} /\partial \mathbf{r}_{j\beta}$ 
can be explicitly represented taking into account that 
\begin{equation*}
\mathbf{f}_{i\alpha} = m_i\mathbf{w}_{i\alpha} - \frac{\partial u(\mathbf{r}_i,t)}{\partial \mathbf{r}_{i \alpha}},
\end{equation*}
 where 
\begin{equation*}
\mathbf{w}_{i\alpha} = -\frac{1}{m_i} \sum_{j(j\neq i)}  \varphi '(r_{ij}) \,\frac{\mathbf{r}_{i \alpha}-\mathbf{r}_{j \alpha})}{r_{ij}}
\end{equation*}
is the inter-particle part of the acceleration.
 The result is 
\begin{equation}\label{eq:fg_term}
 \mathbf{g}_i = -2\sum\limits^{N}_{j(j \neq i)} \Bigl[ (\mathbf{w}_i -\mathbf{w}_j )\frac{\varphi'_{ij}}{r_{ij}} +
\frac{\mathbf{r}_{ij}}{r^{3}_{ij}} \left( r_{ij}\varphi''_{ij} -\varphi'_{ij} \right) 
\bigl(\mathbf{r}_i \cdot (\mathbf{w}_i -\mathbf{w}_j) \bigr) \Bigr]+\mathbf{h}_i,
\end{equation}
where 
\begin{equation*}
\mathbf{h}_i = \frac{2}{m_i}\sum_{\beta} \frac{\partial u}{\partial\mathbf{r}_{i\beta}}\, 
\frac{\partial^{2} u}{\partial\mathbf{r}_{i\alpha}\partial\mathbf{r}_{i\beta}}.
\end{equation*}

Basically the evolution operators $\e^{\A a_{p}h} $ and  
$\e^{\B b_{p}h +\C c_{p}h^3}$ displace $\mathbf{v}$ and $\mathbf{r}$ forward 
in time with  
\begin{equation}\label{eq:shifts}
\begin{aligned}
  &\mathbf{v} \to  \mathbf{v} + b_{p}\mathbf{a} h + c_{p}\mathbf{G}h^{3} & \text{~~~~and }
 & \mathbf{r} \to   \mathbf{r} + a_{p} \mathbf{v} h.
\end{aligned}
\end{equation}
The decomposition integration of Eqn.~\eqref{eq:decomposition} conserves the 
symplectic map of flow of the particles 
in phase space, because the separate shifts of Eqn.~\eqref{eq:shifts} of  
positions and  velocities do not change the phase volume.
Time-reversibility can be ensured by imposing two 
conditions, namely $a_1=0$, $a_{p+1}=a_{P-p+1}$, $b_{p}=b_{P-p+1}$, $c_{p}=c_{P-p+1}$,
as well as $a_{p}=a_{P-p+1}$, $b_{p}=b_{P-p}$, $c_{p}=c_{P-p}$ with $b_P=0$ and $c_P=0$.

Next we deal with numerical integrators of the form given in 
Eqn.~\eqref{eq:decomposition}, the most efficient version of which is due to 
Omelyan \cite{OmMrFo03}.
Adding the force-gradient term $\C$ in the leap-frog scheme does not increase 
the order of the method as one cannot cancel the commutator 
$\bigl[\A,[\A,\B]\bigr]$.  However, the second-order five-stage method 
\begin{equation*}
\Delta(h)_5 = \e^{ \frac{1}{6} h \hat{\B} } \e^{ \frac{1}{2} h \hat{\A} }  \e ^{ \frac{2}{3} h \hat{\B} } \e^{ \frac{1}{2} h \hat{\A} } \e^{ \frac{1}{6} h \hat{\B} },
\end{equation*}
 conserves the shadow Hamiltonian 
\cite{KeSipCl10} 
\begin{equation*}
 \tilde{H}_5 = H -  \bigr[\B,[\A,\B]\bigr]\frac{h^2}{72} + \landO(h^{4}),
\end{equation*}
 where the leading error coefficient is a scalar multiple of the force-gradient term $\C$. Thus  adding a proper amount of the shadow Hamiltonian defines the force-gradient scheme 
\begin{equation*}
\Delta(h)_{5C} = \e^{ \frac{1}{6} h \hat{\B} } \e^{ \frac{1}{2} h \hat{\A} }  \e ^{ \frac{2}{3} h \hat{\B} + \frac{1}{72} h^{3} \C } 
 \e^{\frac{1}{2} h \hat{\A} }  \e^{ \frac{1}{6} h \hat{\B} }. 
\end{equation*}
This scheme conserves the
shadow Hamiltonian \cite{KeSipCl09}
 \begin{equation*}
\begin{split}
 \tilde{H}_{5C} = H +& \left( 41\biggl[\B,\Bigl[\B,\bigl[\B,[\B,\A]\bigr]\Bigr]\biggr] + \ldots \right. \\
   &\left. + 54\biggl[\A,\Bigl[\A,\bigl[\A,[\B,\A]\bigr]\Bigr]\biggr] \right)\frac{h^{4}}{155520}  + \landO(h^{6}),
\end{split}
\end{equation*}
which gains two orders of accuracy.

\subsection{Nested integrators for multirate systems}
In order to reduce the  computational effort to evaluate an evolution operator for one part of the action, we  use a \textit{nested integrator} with a small step-size  
to evaluate the inner cheap part \cite{KeSiCl08}. 
An example of such class of problems can be the multi-time scale problems. 

Let us consider a Hamiltonian which can be represented in the following form  
\begin{equation} \label{eq:multiscale}
 H = T + V_{1}+ V_{2},
\end{equation}
where $T$ represents the kinetic part, $V_{1}$ is the potential energy of the small (fast) scale part of the system and 
$V_{2}$ corresponds to the potential energy of the large (slow) scale part.

We choose the following integrator to compute the inner part  $H = T + V_{1}$ 
\begin{equation*} 
\Delta(h)_{M} = \left[ \e^ {\frac{h}{2M}\hat{\B}_{1}}  \e^{\frac{h}{M}\hat{\A}}  \e^{\frac{h}{2M}\hat{\B}_{1}} \right]^  {M}. 
\end{equation*}

Therefore  we  define $ \Delta\left(h\right)$  a nested integrator to solve the  split problem  of Eqn.~\eqref{eq:multiscale},
it yields
\begin{equation}\label{eq:integrator} 
\hat{\Delta}(h) = \left[ \e^{\frac{h}{2} \hat{\B}_{2}}  \Delta\left(h\right)_{M} \e^{\frac{h}{2}\hat{\B}_{2}} \right]^{l}.
\end{equation} 

This method, called \textit{nested leap-frog}, conserves the  shadow Hamiltonian  \cite{KeSiCl08}
 \begin{align*}
 \tilde{H}_M =  H +  \biggl( 
                                -\frac{1}{24}\bigl[\B_2,[\B_2,\A]\bigr] + \frac{1}{12}\bigl[\B_1,[\B_2,\A]\bigr] 
                                +\frac{1}{12}\bigl[\A,[\B_2,\A]\bigr] \\ 
             + \frac{1}{M^{2}}  \Bigl( -\frac{1}{24}\bigl[\B_1,[\B_1,\A]\bigr] +\frac{1}{12}\bigl[\A,[\B_1,\A]\bigr] \Bigr) \biggr) h^{2} + \landO(h^{4}). 
  \end{align*}

\section{Combining Force-Gradient and Multirate Splitting Technique}

Our idea is to combine both the  force-gradient  and the nested algorithm approaches in order to obtain a higher 
energy conservation rate. To do so, let us first  take a look at the 
 following \textit{alike 5-stage nested  integrator}
\begin{equation} \label{eq:multi_force1} 
\Delta(h) =  \left[ \e^{ \lambda h \hat{\B}_{2}}  \Delta\left( \frac{h}{2} \right)_{M}  \e^{  (1- 2\lambda ) h \hat{\B}_{2}}   
\Delta\left(\frac{h}{2}\right)_{M} \e^{\lambda h\hat{\B}_{2} } \right]^{l},
\end{equation}
where 
\begin{equation*}  
\Delta\left(\frac{h}{2}\right)_M =  \left[ \e^{\frac{h}{4M} \hat{\B}_{1}}  \e^{\frac{h}{2M}\hat{\A}} \e^{ \frac{h}{4M}\hat{\B}_{1}} \right]^{M}.
\end{equation*}
We have chosen the 5-stage numerical integrator, 
since it has an optimal number of steps, necessary for  increasing its order.
To analyze the energy conservation of this integrator we have to determine its shadow Hamiltonian. 

In order to do so, we use the BCH formula \eqref{eq:BCH}. 
To simplify this task we consider the limit of the integrator of 
Eqn.~\eqref{eq:multi_force1}, as $M$ tends to infinity. We obtain 
\begin{equation} \label{eq:multi_force_lim} 
\Delta(h) =  \left[ \e^ {\lambda h \hat{\B}_{2}}  \e^{\frac{h}{2}(\hat{\B}_1 +\A) }  \e^{ (1- 2\lambda ) h \hat{\B}_{2} }   
\e^{\frac{h}{2}(\hat{\B}_1 +\hat{\A}) }  \e^{\lambda h\hat{\B}_{2}} \right]^{l}.
\end{equation}

\begin{theor}[Shadow Hamiltonian of \eqref{eq:multi_force_lim}]
The shadow Hamiltonian of the nested multirate integrator \eqref{eq:multi_force_lim} is given by
\begin{multline}
 \tilde{H} =  H  + 
\Bigl( \frac{-1+6\lambda-6\lambda^{2}}{12} \bigl[\B_2,[\A,\B_2]\bigr] \\
+\frac{-1+6\lambda}{24} \bigl[\B_1,[\A,\B_2]\bigr]      \frac{-1+6\lambda}{24} \bigl[\A,[\A,\B_2]\bigr]  \Bigr) h^{2} + \landO(h^{4}). \label{shad.ham.ourmeth}
\end{multline}
\end{theor}
{\em Proof:} 
We apply the BCH formula to the first two evolution operators 
 \begin{equation*}
X=\ln \left( \e^{\lambda h \hat{\B}_2}\e^{ h\frac{ \hat{\B}_1 +\hat{\A}}{2}} \right)= c_1 h +c_2 h^{2}+ c_3 h^{3} + \landO(h^{5}),
\end{equation*}
where
\begin{align*}
  \mathbf{c_1} &= \lambda  \B_2 + \frac{ \B_1 +\A}{2}, \\
  2 c_2       &= \frac{B_2}{2!} \text{ad}c_1 \left(\lambda  \B_2 + \frac{ \B_1 +\A}{2}\right) - \frac{1}{2} \text{ad} c_1\left(\lambda  \B_2 - \frac{ \B_1 +\A}{2}\right),\\
  \mathbf{c_2}&= - \frac{1}{4}\left[\lambda  \B_2 + \frac{ \B_1 +\A}{2},~ \lambda  \B_2 - \frac{ \B_1 +\A}{2} \right] = - \frac{\lambda}{4}[\B_1,\B_2]- \frac{\lambda}{4}[\A,\B_2],\\
  3 c_3       &= \frac{B_3}{3!} \text{ad}c_1 \text{ad}c_1  \left(\lambda  \B_2 + \frac{ \B_1 +\A}{2}\right) - \frac{1}{2} \text{ad} c_2\left(\lambda  \B_2 - \frac{ \B_1 +\A}{2}\right),\\
 \mathbf{c_3} &= - \frac{1}{6}\left[- \frac{\lambda}{4}[\B_1,\B_2]- \frac{\lambda}{4}[\A,\B_2],~ \lambda  \B_2 - \frac{ \B_1 +\A}{2} \right] \\
	      &= 
- \frac{\lambda^{2}}{24} \bigl[\B_2,[\B_1,\B_2]\bigr] 
 - \frac{\lambda^{2}}{24} \bigl[\B_2,[\A,\B_2]\bigr]
+ \frac{\lambda}{48} \bigl[\B_1,[\B_1,\B_2]\bigr]\\ 
              & \qquad
+ \frac{\lambda}{48} \bigl[\B_1,[\A,\B_2]\bigr] 
+ \frac{\lambda}{48} \bigl[\A,[\B_1,\B_2]\bigr] 
+ \frac{\lambda}{48} \bigl[\A,[\A,\B_2]\bigr].
\end{align*}

Then we have the result for our first two operators 
\begin{align*}
 X &= \left(\lambda  \B_2 + \frac{ \B_1 +\A}{2} \right) h 
+ \left(- \frac{\lambda}{4}[\B_1,\B_2]- \frac{\lambda}{4}[\A,\B_2]\right) h^{2} \\
   &\qquad
+\left(- \frac{\lambda^{2}}{24} \bigl[\B_2,[\B_1,\B_2]\bigr] 
 - \frac{\lambda^{2}}{24} \bigl[\A,[\A,\B_2]\bigr] 
 + \frac{\lambda}{48} \bigl[\B_1,[\B_1,\B_2]\bigr] \right. \\
   &\left. \qquad 
+\frac{\lambda}{48} \bigl[\B_1,[\A,\B_2]\bigr] 
+ \frac{\lambda}{48} \bigl[\A,[\B_1,\B_2]\bigr] 
 + \frac{\lambda}{48} \bigl[\A,[\A,\B_2]\bigr]  \right) h^{3} + \landO(h^{5}).
\end{align*}

The next step is to apply the BCH formula on the following operators
\begin{equation*}  
Y= \ln \left( \e^{X}\e^{  (1- 2\lambda ) h \hat{\B}_{2}}  \right)= c_1 +c_2 + c_3  + \landO(h^{5})
\end{equation*}
and coefficients
\begin{align*}
 \mathbf{c_1} &=X + (1- 2\lambda ) h \B_{2} = \left((1-\lambda )  \B_{2} +\frac{\B_1+\A}{2}\right)h + (\ldots)h^{2} +(\ldots)h^{3}, \\
  2 c_2 &= \frac{B_2}{2!} \text{ad}c_1 \bigl(X +(1- 2\lambda ) h \B_{2}\bigr)
            - \frac{1}{2} \text{ad} c_1\bigl(X -(1- 2\lambda ) h \B_{2}\bigr),\\
  \mathbf{c_2} &= - \frac{1}{4}\bigl[X +(1- 2\lambda ) h \B_{2} ,~ X -(1- 2\lambda ) h \B_{2}\bigr] \\ 
               &=\left(\frac{1-2\lambda}{4} [\B_1,\B_2] +\frac{1-2\lambda}{4} [\A,~\B_2]\right) h^{2} \\
               &\qquad+ \left(\frac{\lambda-2\lambda^{2}}{8}\bigl[\B_2,[\B_1,\B_2]\bigr]
                      +\frac{\lambda-2\lambda^{2}}{8}\bigl[\B_2,[\A,\B_2]\bigr]\right) h^{3}\\
    3 c_3      &= \frac{B_3}{3!} \text{ad}c_1 \text{ad}c_1  \bigl(X +(1- 2\lambda ) h \B_{2}\bigr)
                   - \frac{1}{2} \text{ad}c_2\bigl(X -(1- 2\lambda ) h \B_{2}\bigr),\\
 \mathbf{c_3} &= - \frac{1}{6}\bigl[c_2,~ X -(1- 2\lambda ) h \B_{2} \bigr]\\
& =  \left( 
- \frac{(1- 2\lambda )(1- 3\lambda )}{24} \bigl[\B_2,[\B_1,\B_2]\bigr] 
- \frac{(1- 2\lambda )(1- 3\lambda )}{24} \bigl[\B_2,[\A,\B_2]\bigr]\right. \\
 &\left.\qquad  + \frac{(1- 2\lambda )}{48} \bigl[\B_1,[\B_1,\B_2]\bigr] 
+ \frac{(1- 2\lambda )}{48} \bigl[\B_1,[\A,\B_2]\bigr]
 + \frac{(1- 2\lambda )}{48} \bigl[\A,[\B_1,\B_2]\bigr] \right. \\
 &\left.\qquad+ \frac{(1- 2\lambda )}{48} \bigl[\A,[\A,\B_2]\bigr] \right) h^{3},
\end{align*}
and we obtain the following expansion 

\begin{align*}
 Y&= \left((1-\lambda ) \B_2 + \frac{ \B_1 +\A}{2} \right) h +\left( \frac{1-3\lambda}{4}[\B_1,\B_2]+\frac{1-3\lambda}{4}[\A,\B_2]\right) h^{2} \\
& \qquad 
+\left(\frac{-1+8\lambda-13\lambda^{2}}{24} \bigl[\B_2,[\B_1,\B_2]\bigr] 
     + \frac{-1+8\lambda-13\lambda^{2}}{24} \bigl[\A,[\A,\B_2]\bigr]   \right. \\
 &\left. \qquad
  + \frac{1-\lambda}{48} \bigl[\B_1,[\B_1,\B_2]\bigr]
  + \frac{1-\lambda}{48} \bigl[\B_1,[\A,\B_2]\bigr] 
  + \frac{1-\lambda}{48} \bigl[\A,[\B_1,\B_2]\bigr] \right. \\
 &\left. \qquad
+ \frac{1-\lambda}{48} \bigl[\A,[\A,\B_2]\bigr]  \right) h^{3} + \landO(h^{5}),
\end{align*}

The next step would be to repeat the previous procedures to find
\begin{equation*}
 Z = \ln \left( \e^{Y}\e^{  h\frac{ \hat{\B}_1 +\hat{\A}}{2} }  \right)= c_1 +c_2 + c_3  + \landO(h^{5}).
\end{equation*}

Using the BCH formula we obtain 
\begin{align*}
 \mathbf{c_1} &=Y + h\frac{ \B_1 +\A}{2} = ((1-\lambda )  \B_{2} + \B_1+\A)h + (\ldots)h^{2} +(\ldots)h^{3}, \\
  2 c_2 &= \frac{B_2}{2!} \text{ad}c_1 \left(Y +h\frac{ \B_1 +\A}{2} \right) - \frac{1}{2} \text{ad} c_1\left(Y - h\frac{ \B_1 +\A}{2}\right),\\
  \mathbf{c_2} &= - \frac{1}{4}\left[Y +h\frac{ \B_1 +\A}{2} ,~ Y -h\frac{ \B_1 +\A}{2}\right] \\ 
 &=\left(-\frac{1-\lambda}{4} [\B_1,~\B_2] -\frac{1-\lambda}{4} [\A,~\B_2]\right) h^{2}\\
  & \qquad 
+\left(-\frac{1-3\lambda}{16}\bigl[\B_1,[\B_1,\B_2]\bigr] 
       -\frac{1-3\lambda}{16}\bigl[\B_1,[\A,\B_2]\bigr] \right. \\
&\left. \qquad \qquad-\frac{1-3\lambda}{16}\bigl[\A,[\B_1,\B_2]\bigr] 
               -\frac{1-3\lambda}{16}\bigl[\A,[\A,\B_2]\bigr]\right) h^{3}, \\
  3 c_3 &= \frac{B_3}{3!} \text{ad}c_1 \text{ad}c_1  \left(Y +h\frac{ \B_1 +\A}{2}\right) - \frac{1}{2} \text{ad} c_2\left(Y -h\frac{ \B_1 +\A}{2}\right)\\
 \mathbf{c_3} &= - \frac{1}{6}\left[c_2,~ Y - h\frac{ \B_1 +\A}{2} \right],\\ 
&=  \left( - \frac{(1- \lambda )^{2}}{24} \bigl[\B_2,[\B_1,\B_2]\bigr] 
          - \frac{(1- \lambda )^{2}}{24} \bigl[\B_2,[\A,\B_2]\bigr] \right) h^{3}.
\end{align*}
Therefore we obtain 
\small
\begin{align*}
 Z&= \bigl((1-\lambda ) \B_2 + \B_1 +\A \bigr) h 
   +\Bigl( \frac{-\lambda}{2}[\B_1,\B_2]+\frac{-2\lambda}{2}[\A,\B_2]\Bigr) h^{2} \\
  &\qquad 
+ \left(\frac{-1+5\lambda-7\lambda^{2}}{12} \bigl[\B_2,[\B_1,\B_2]\bigr]
      + \frac{-1+5\lambda-7\lambda^{2}}{12} \bigl[\A,[\A,\B_2]\bigr]\right.\\
  &\left.\qquad\qquad
+ \frac{-1+4\lambda}{48} \bigl[\B_1,[\B_1,\B_2]\bigr] 
+ \frac{-1+4\lambda}{48} \bigl[\B_1,[\A,\B_2]\bigr] \right. \\
  &\left.\qquad\qquad
+ \frac{-1+4\lambda}{48} \bigl[\A,[\B_1,\B_2]\bigr] 
  +\frac{-1+4\lambda}{48} \bigl[\A,[\A,\B_2]\bigr]  \right) h^{3} + O(h^{4}).
\end{align*}
\normalsize
Applying the BCH formula for a last time we obtain the shadow Hamiltonian 
\begin{equation*}
\tilde{H}  = \ln \left( \e^{Z}\e^{  \lambda \hat{\B}_2 h  }  \right)= c_1 +c_2 + c_3  + \landO(h^{5}),
\end{equation*}
with the coefficients\small
\begin{align*}
 \mathbf{c_1} &=Z + \lambda \B_2 h = ( \B_{2} + \B_1+\A)h + (\ldots)h^{2} +(\ldots)h^{3}, \\
  2 c_2 &= \frac{B_2}{2!} \text{ad}c_1 \left(Z + \lambda \B_2 h  \right) - \frac{1}{2} \text{ad} c_1\left(Z - \lambda \B_2 h \right),\\
  \mathbf{c_2} &= - \frac{1}{4}\left[Z + \lambda \B_2 h  ,~ Z - \lambda \B_2 h \right] \\
&= 
 \left(\frac{\lambda}{2} [\B_1,~\B_2] +\frac{\lambda}{2} [\A,~\B_2]\right) h^{2} 
 +\left(\frac{\lambda^2}{4}\bigl[\B_2,[\B_1,\B_2]\bigr] 
   +\frac{\lambda^2}{4}\bigl[\B_2,[\A,\B_2]\bigr] \right) h^{3},\\
  3 c_3 &= \frac{B_3}{3!} \text{ad}c_1 \text{ad}c_1  \left(Z + \lambda \B_2 h\right) - \frac{1}{2} \text{ad} c_2\left(Z - \lambda \B_2 h\right),\\
 \mathbf{c_3} &= - \frac{1}{6}\left[c_2,~ Z - \lambda \B_2 h \right] \\
&=  \left(  \frac{\lambda(1- 2\lambda )^{2}}{12} \bigl[\B_2,[\B_1,\B_2]\bigr] 
          + \frac{\lambda(1- 2\lambda)^{2}}{12} \bigl[\B_2,[\A,\B_2]\bigr] 
            +\frac{\lambda }{12} \bigl[\B_1,[\B_1,\B_2]\bigr]  \right.\\
 &\left.  \qquad+\frac{\lambda }{12} \bigl[\B_1,[\A,\B_2]\bigr]
  +\frac{\lambda }{12} \bigl[\A,[\B_1,\B_2]\bigr]
  +\frac{\lambda }{12} \bigl[\A,[\A,\B_2]\bigr] \right) h^{3},
\end{align*}
\normalsize
and finally
\small
\begin{align*}
 \tilde{H} &= H
 + \left(\frac{-1+6\lambda-6\lambda^{2}}{12} \bigl[\B_2,[\B_1,\B_2]\bigr]+
    \frac{-1+6\lambda}{24} \bigl[\B_1,[\A,\B_2]\bigr] \right.\\
  &\left.\qquad 
+ \frac{-1+6\lambda-6\lambda^{2}}{12} \bigl[\B_2,[\A,\B_2]\bigr] 
 + \frac{-1+6\lambda}{24} \bigl[\B_1,[\B_1,\B_2]\bigr]\right.\\
 &\left. \qquad
+ \frac{-1+6\lambda}{24} \bigl[\A,[\B_1,\B_2]\bigr] 
   +\frac{-1+6\lambda}{24} \bigl[\A,[\A,\B_2]\bigr]  \right) h^{2}+ \landO(h^{4}).
\end{align*}
\normalsize
Finally, taking into account that $[\B_1,\B_2] = 0$, 
hence $\bigl[\B_2,[\B_1,\B_2]\bigr]$, $\bigl[\B_1, [\B_1 , \B_2]\bigr]$ and 
$\bigl[\A, [\B_1 , \B_2]\bigr]$ are equal to zero, we obtain \eqref{shad.ham.ourmeth}. \qed

We can eliminate a couple of  terms by choosing  $\lambda =1/6$, thus 
\begin{equation*}
 \tilde{H} = H    -\frac{1}{72} \bigl[\B_2,[\A,\B_2]\bigr]  h^{2} 
+ \landO(h^{4}).
\end{equation*}
We would like to increase the order of the method \eqref{eq:multi_force1} by adding the force-gradient term,
but  first we consider the force-gradient itself. 
Due to the  splitting \eqref{eq:multiscale}  it can be represented as
\begin{multline*}
 C = \bigl[\B,[\A,\B]\bigr] = \bigl[\B_2 +\B_1,[\A,\B_2+\B_1]\bigr] \\
= \bigl[\B_2,[\A,\B_2]\bigr] + \bigl[\B_1,[\A,\B_1]\bigr] + \bigl[\B_1,[\A,\B_2]\bigr] + \bigl[\B_2,[\A,\B_1]\bigr].
\end{multline*}
\normalsize
Then we can tune the original algorithm \eqref{eq:multi_force1} by 
adding the first term of  the force gradient $\bigl[\B_2,[\A,\B_2]\bigr]$ and neglect
 the last three terms:
\small
\begin{equation} \label{eq:multi_force2} 
\Delta(h) =  \left[ \e ^{\frac{1}{6} h \hat{\B}_{2}} \Delta\left( \frac{h}{2} \right)_{M}  \e^{ \frac{2}{3} h \hat{\B}_{2} + 
\frac{1}{72}h^{3}\bigl[\B_2,[\A,\B_2]\bigr]}    \Delta\left(\frac{h}{2}\right)_{M}  \e^{\frac{1}{6} h\hat{\B}_{2}} \right]^{l},
\end{equation}
\normalsize
which preserves the fourth-order accurate shadow Hamiltonian 
\begin{equation*}
 \tilde{H} = H    + \landO(h^{4}).
\end{equation*}

\section{Numerical Experiments}
In order to estimate the performance of  the integrator of Eqn.~\eqref{eq:multi_force2} 
we compare it with the other algorithms mentioned above. 
Let us consider the three body problem \cite{Hairer} and a particular case of it,
the   \textit{Sun-Earth-Moon problem}. 
The given system has the energy 
\begin{equation*}
 E=\sum\limits^{2}_{i=0} \frac{m_{i}v^2_{i}}{2} - G \sum\limits^{2}_{i=1}\sum\limits^{i-1}_{j=0} \frac{m_{i}m_{j}}{r_{ij}}, 
\end{equation*}
 where $r_{ij}=\|\mathbf{r_{i}-r_{j}}\|$, $m_{0}$, $m_{1}$ and $m_{2}$ represent the masses of the Sun, the Earth
 and the Moon, respectively and $G$ is the gravitational constant.
The equations of motion are then
\begin{equation} \label{eq:three_body}
\begin{aligned}
 &\frac{\dd \mathbf{r}_{0} }{\dd t} = \mathbf{v}_{0}, \qquad
&\frac{\dd\mathbf{v}_{0}}{\dd t}=-  m_{1}G\,\frac{\mathbf{r}_{0}-\mathbf{r}_{1}}{r_{01}^3}
                                  - m_{2}G\,\frac{\mathbf{r}_{0}-\mathbf{r}_{2}}{r_{02}^3},\\
 &\frac{\dd \mathbf{r}_{1} }{\dd t} = \mathbf{v}_{1}, &\frac{\dd\mathbf{v}_{1}}{\dd t}
=-m_{0}G\,\frac{\mathbf{r}_{1}-\mathbf{r}_{0}}{r_{10}^3} 
- m_{2}G\,\frac{\mathbf{r}_{1}-\mathbf{r}_{2}}{r_{12}^3},\\
 &\frac{\dd \mathbf{r}_{2} }{\dd t} = \mathbf{v}_{2}, &\frac{\dd\mathbf{v}_{2}}{\dd t}
 =-m_{0}G\,\frac{\mathbf{r}_{2}-\mathbf{r}_{0}}{r_{20}^3} 
 - m_{1}G\,\frac{\mathbf{r}_{2}-\mathbf{r}_{1}}{r_{21}^3}.
\end{aligned}
\end{equation}
The force-gradient terms can be obtained from  \eqref{eq:fg_term} for  this case, 
using the external field potential $u(r_{ij})=0$ and 
the pair-wise potentials
\begin{equation*}
 \varphi(r_{ij})=-G \frac{m_{i}m_{j}}{r_{ij}},
\end{equation*}
respectively for each interaction.
\begin{table}
{\footnotesize
\begin{tabular}{||l|c|c||}
\hline
  Gravitational constant        $(G)$            & $6.67384 \times 10^{-11}$,$\text{m}^{3}$/kg s & $0.2662$ $\text{AU}^{3}$/SU mo \\ \hline
  Mass of the Sun               $(m_0)$          & $ 1.9891 \times 10^{30}$, kg                & $1$ SU \\ \hline
  Mass of the Earth             $(m_1)$          & $ 5.9736 \times 10^{24}$, kg                & $3 \times 10^{-6}$ SU \\ \hline
  Mass of the Moon              $(m_2)$          & $ 7.3477 \times 10^{22}$, kg                & $0.0369 \times 10^{-6}$ SU \\  \hline
  Initial position of the Sun   $(\mathbf{r}_0)$ & $(0,~0)$,        m                          & $(0,~0)$, AU   \\ \hline
  Initial position of the Earth $(\mathbf{r}_1)$ & $(0,~1.52098 \times 10^{11})$, m            & $(0,~1.0167138)$, AU \\ \hline
  Initial position of the Moon  $(\mathbf{r}_2)$ & $(0,~1.52504 \times 10^{11})$, m            & $(0,~1.0191138)$, AU \\ \hline
  Initial velocity of the Sun   $(\mathbf{v}_0)$ & $(0,~0)$,                      m/s          & $(0,~0)$, AU/mo \\ \hline
  Initial velocity of the Earth $(\mathbf{v}_1)$ & $(0,~29.78 \times 10^{3})$,    m/s          & $(0,~0.5160)$, AU/mo \\ \hline
  Initial velocity of the Moon  $(\mathbf{v}_2)$ & $(0,~ 30.802\times 10^{3})$,   m/s          & $(0,~0.5337)$, AU/mo \\ \hline
\end{tabular} 
}
\caption{Physical parameters of the Sun-Earth-Moon problem.}
\end{table}

\begin{figure}[h!]
\centering
\includegraphics[width=0.49\linewidth]{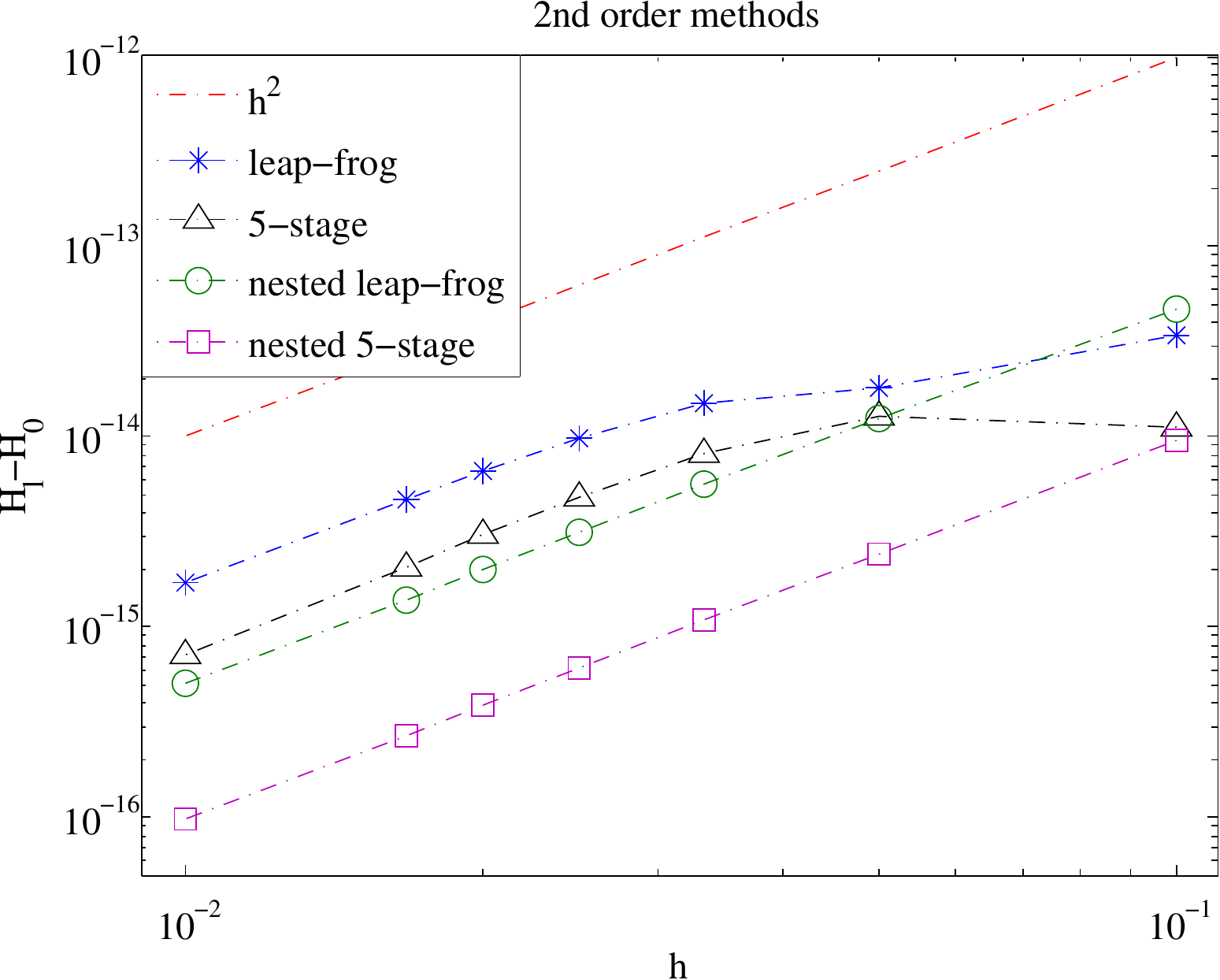} 
\includegraphics[width=0.49\linewidth]{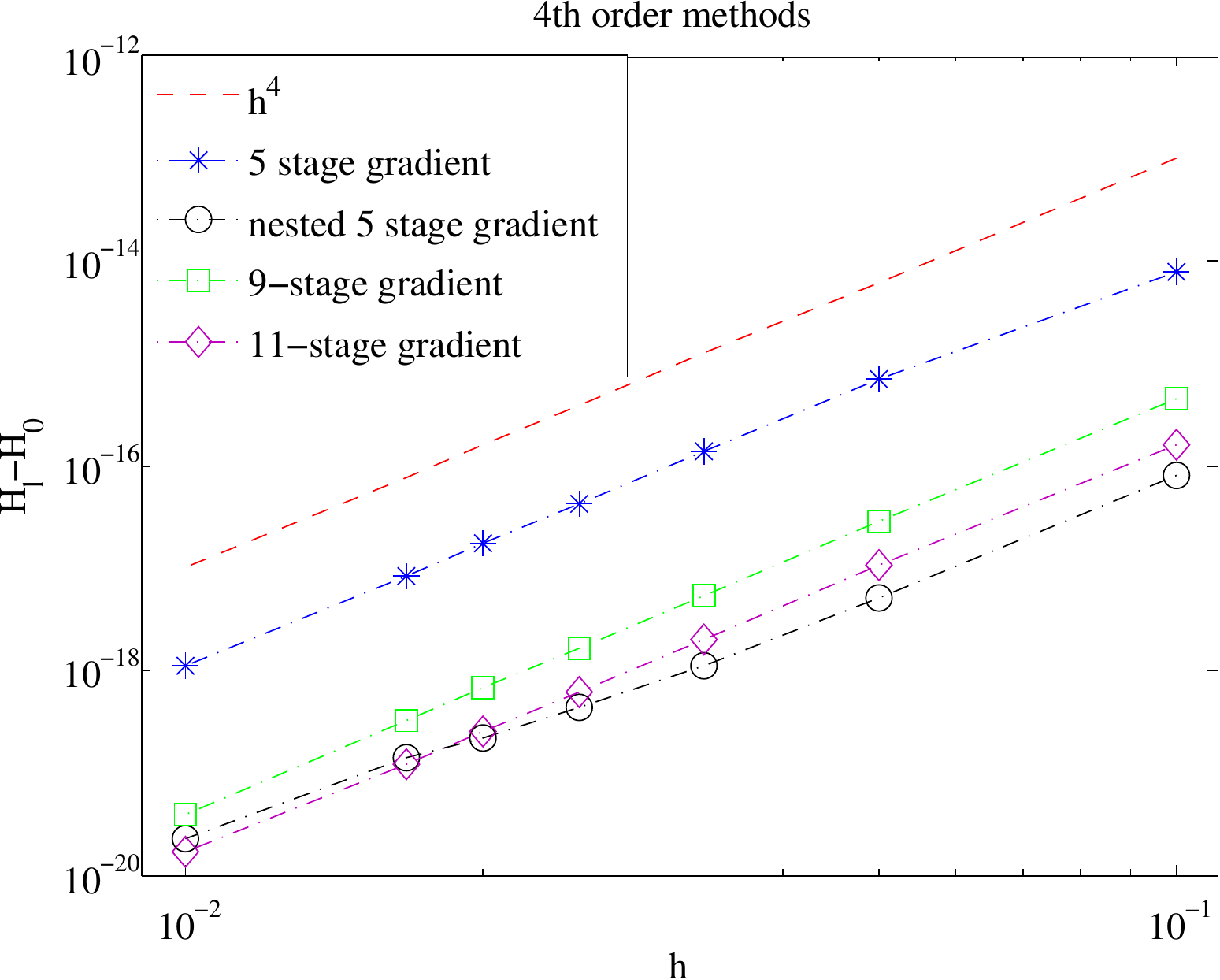}  \vspace{-10pt}
  \caption{Sun-Moon-Earth problem: absolute error for different integrators}
\label{fig:1}
\end{figure}

Figure~\ref{fig:1} presents a comparison between the standard numerical 
algorithms, nested approaches, the force-gradient and our combined method. 
The proposed integrator of Eqn.~\eqref{eq:multi_force2} with $M=30$, which 
combines nested and force-gradient ideas, yields a better energy conservation 
even compared with 9-stage and 11-stage force-gradient numerical schemes.  
These numerical results correspond to our analytical observations.

\begin{figure}[!ht]
\centering
\includegraphics[width=1.0 \linewidth]{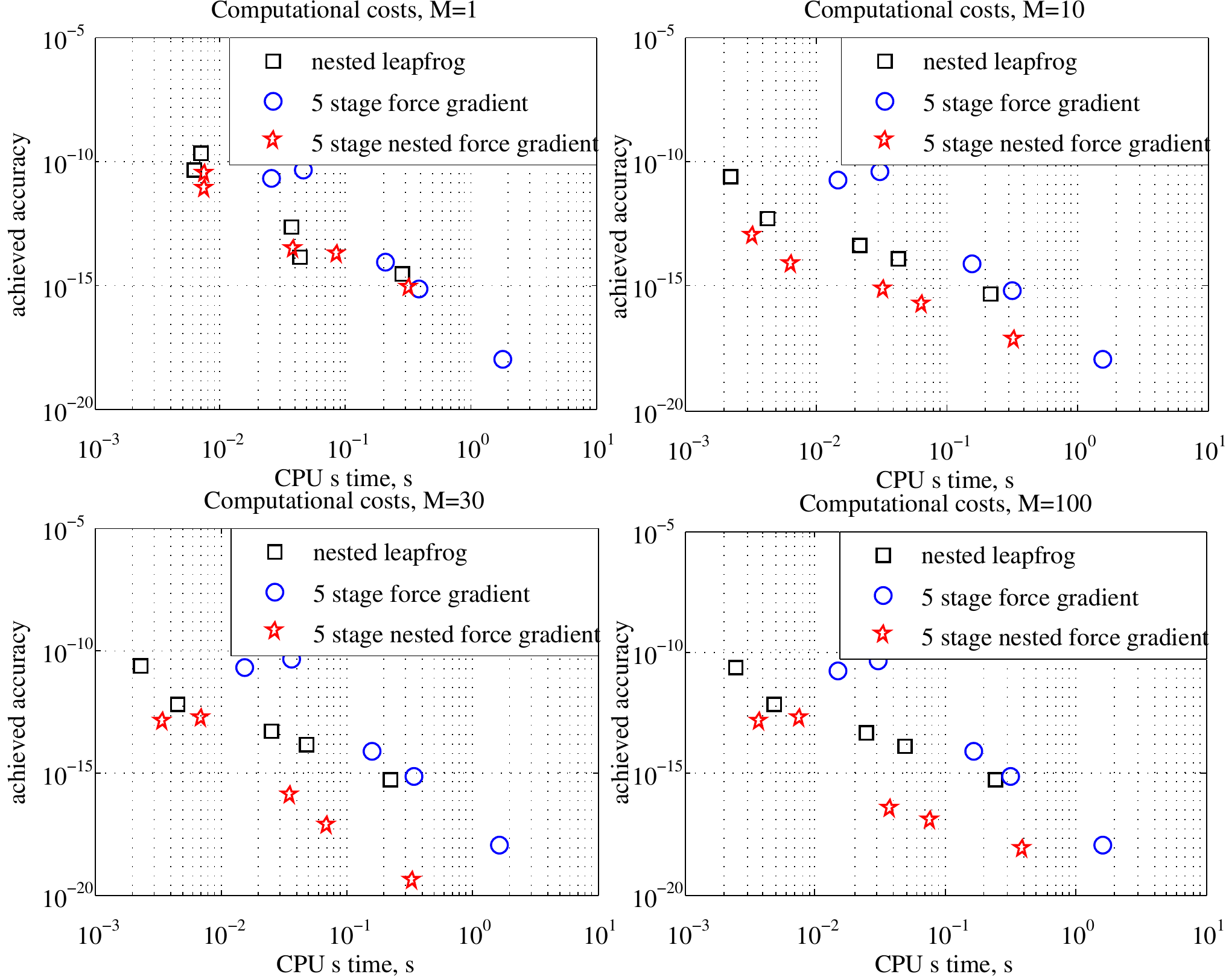}  \vspace{-10pt}
  \caption{Sun-Moon-Earth problem: CPUs time  vs. achieved accuracy for different integrators}
\label{fig:2}
\end{figure}

Figure~\ref{fig:2} presents the CPU time, required for the three different 
integrators against the achieved 
accuracy. Here we scale the time needed for the computation of the fast part by a factor of $0.001$, since we assume that 
the computation of the fast scale functions is very cheap compared to the slow 
scale function evaluations. 
We can see that in general our nested force-gradient method \eqref{eq:multi_force2} requires less CPU time and performs more
 accurate than the standard schemes, presented in Figure~\ref{fig:2}.\\

Thus we can argue that, if the evaluation of fast 
function is significantly cheaper than the {\it slow} function, computational costs decrease. This is exactly the case found in 
our long-term goal applications in lattice quantum chromodynamics (LQCD), 
where the action can be split into two parts: the gauge 
action~(whose force evaluations are cheap) and the fermion action~(expensive).  
  

\section{Conclusions and Outlook}
We have introduced a new decomposition scheme for Hamiltonian systems, 
which combines the idea of the force-gradient time-reversible and 
symplectic integrators and 
the splitting approach of nested algorithms. The new
 method of Eqn.~\eqref{eq:multi_force2} is fourth-order accurate. Compared to other fourth-order schemes, the leading error coefficient is smaller and computational costs are lower.

Our future work will apply this approach in the Hybrid Monte Carlo \cite{HMC} 
(HMC) 
algorithm for numerical integration of the lattice path-integral of 
quantum chromodynamics (QCD), which describes the strong interactions between 
quarks and gluons inside the nucleons. In this case, the Hamiltonian dynamics 
are defined on curved manifolds and one has to take into account the 
non-commutativity of the operators $\A$ and $\B$.


\section*{Acknowledgments}
{
This work is supported by the European Union within the Marie Curie Initial Training Network STRONGnet
 on {\em  Strong Interaction Supercomputing Training Network} (Grant Agreement number 238353).}
This work is part of project B5 within the SFB/Transregio 55
{\em Hadronenphysik mit Gitter-QCD}.


\end{document}